\documentclass{amsart}

\def\EE{{\mathcal E}}
\def\E{{\mathbb E}}
\def\var{{\rm Var}}
\def\ent{{\rm Ent}}
\def\R{{\sf R}}

\theoremstyle{plain}
\newtheorem{lemma}{Lemma}[section]
\newtheorem{corollary}{Corollary}[section]
\newtheorem{theorem}{Theorem}[section]

\theoremstyle{definition}
\newtheorem{definition}{Definition}[section]

\newtheorem{remark}{Remark}[section]

\numberwithin{equation}{section}

\begin{document}

\title{Mixing Time bounds via the Spectral Profile}
\author{Sharad Goel, Ravi Montenegro and Prasad Tetali}
\thanks{Research supported in part by NSF grants DMS-0306194, DMS-0401239}
\subjclass[2000]{60,68}


\begin{abstract}
On complete, non-compact manifolds and infinite graphs, Faber-Krahn inequalities have been used to estimate the rate of decay of the heat kernel. We develop this technique in the setting of finite Markov chains, proving upper and lower $L^{\infty}$ mixing time bounds via the spectral profile. This approach lets us recover and refine previous conductance-based bounds of mixing time (including the Morris-Peres result), and in general leads to sharper estimates of convergence rates. We apply this method to several models including groups with moderate growth,  the fractal-like Viscek graphs, and the product group
$\mathbb{Z}_a \times \mathbb{Z}_b$, to obtain tight bounds on the corresponding mixing times.
\end{abstract}

\maketitle

\section{Introduction} \label{sec:intro}

It is well known that the spectral gap of a Markov chain can be estimated in terms of conductance, facilitating isoperimetric bounds on mixing time (see \cite{JSiso, LawSok}). Observing that small sets often have large conductance, Lov{\'a}sz and Kannan (\cite{LovCon}) strengthened this result by bounding total variation mixing time for reversible chains in terms of the ``average conductance'' taken over sets of various sizes. Morris and Peres (\cite{Morris}) introduced the idea of evolving sets to analyze reversible and non-reversible chains, and found bounds on the larger $L^{\infty}$ mixing time.

To sidestep conductance, we introduce ``spectral profile" and develop Faber-Krahn inequalities in the context of finite Markov chains, bounding mixing time directly in terms of the spectral profile. FK-inequalities were introduced by Grigor'yan and developed together with Coulhon and Pittet (\cite{Grig1, Coulhon, Grig2, Barlow}) to estimate the rate of decay of the heat kernel on manifolds and infinite graphs. Their techniques build on functional analytic methods presented, for example, in \cite{Davies}. We adapt this approach to the setting of finite Markov chains and derive $L^{\infty}$ mixing time estimates for both reversible and non-reversible walks.

These bounds let us recover the previous conductance-based results, and in general lead to sharper estimates on rates of convergence to stationarity. We also show that the spectral profile can be bounded in terms of both log-Sobolev and Nash inequalities, leading to new and elementary proofs for previous mixing time results -- for example, we re-derive Theorem 3.7 of Diaconis--Saloff-Coste \cite{LogSob} and   Theorem 42 (Chapter 8) of Aldous-Fill \cite{AldousBook}.

In terms of applications, we  first observe that for simple examples such as  the random walk on a complete graph and the $n$-cycle, the spectral profile gives the correct bounds.  As more interesting examples, we also analyze walks on graphs with moderate growth,  the fractal-like Viscek graphs, product groups like $Z_a \times Z_b$, and show optimal bounds. In the case of the graphs with moderate growth, we show that the mixing time is of the order of the square of the diameter,  a result originally due to Diaconis and Saloff-Coste (see \cite{ModGrowth, DSNash}).  
In the case of the Viscek graphs, we show that the spectral profile provides tight upper and lower bounds on mixing time, and observe that the conductance-based bounds give much weaker upper bounds. 

In Section~\ref{sec:intro} we introduce notation, review preliminary ideas and state our main results. Section~\ref{sec:bounds} presents the proofs of both the continuous and discrete time versions of the spectral profile upper bound on mixing time. In Section~\ref{sec:lower} we recall a complementary lower bound shown in \cite{Grig2}. Section~\ref{sec:apps} discusses applications, including the relationship between the spectral profile, log-Sobolev and Nash inequalities. Section~\ref{sec:viscek} discusses the more elaborate example of the Viscek graphs. Section~\ref{sec:delicate_eg} discusses the spectral profile of the random walk on $Z_a \times Z_b$, which turns out to be  a bit subtle.

\subsection{Preliminaries}
A Markov chain on a finite state space $\mathcal{X}$ can be identified with a kernel $K$ satisfying
$$ K(x,y) \geq 0 \hspace{.5cm} \sum_{y \in \mathcal{X}} K(x,y) = 1. $$
The kernel of $K^n$ is then given iteratively by 
$\displaystyle K_n(x,y) = \sum_{z \in \mathcal{X}} K_{n-1}(x,z) K(z,y)$
and can be interpreted as the probability of moving from state $x$ to $y$ in exactly $n$ steps. 

We say that a probability measure $\pi$ on $\mathcal{X}$ is invariant with respect to $K$ if 
$\displaystyle  \sum_{x \in \mathcal{X}} \pi(x) K(x,y) = \pi(y). $
That is, starting with distribution $\pi$ and moving according to the kernel $K$ leaves the distribution of the chain unchanged. Throughout, we assume that $K$ is irreducible: For each $x,y \in \mathcal{X}$ there is an $n$ such that $K_n(x,y) > 0$. Under this assumption $K$ has a unique invariant measure $\pi$ and $\pi_* = \min_{x} \pi(x) > 0$. 

The chain $(K, \pi)$ is reversible if, 
$K=K^*$ is a self-adjoint operator on the Hilbert space $L^2(\pi)$.
In general, $K^*(x,y)=\pi(y)K(y,x)/\pi(x)$, and so reversibility
is equivalent to requiring that $K$ satisfy the detailed balance equation:
for all $x,y \in \mathcal{X}$, we have  $\pi(x) K(x,y) = \pi(y) K(y,x)$.

The kernel $K$ describes a discrete-time chain which at each step moves with distribution according to $K$. Alternatively, we can consider the continuous-time chain $H_t$, which waits an exponential time before moving. More precisely, as operators
$$H_t = e^{-t\Delta}  \hspace{.5cm} \Delta = I-K.$$ 
The kernel of $H_t$ is then given explicitly by
$$H_t(x,y) = e^{-t}\sum_{n=0}^{\infty} \frac{t^n}{n!}K^n(x,y).$$
Let $h(x,y,t) = H_t(x,y)/\pi(y)$ denote the density of $H_t(x, \cdot)$ with respect to its stationary measure $\pi$. 

%

To measure the rate of convergence to equilibrium, we first need to decide on a distance.
\begin{definition}
For two measures $\mu$ and $\nu$ with densities $f(x) = \mu(x)/\pi(x)$ and $g(x) = \nu(x)/\pi(x)$ with respect to the positive measure $\pi$, their $L^p(\pi)$ distance is 
$$d_{p, \pi}(\mu, \nu) = \|f - g\|_{L^p(\pi)} \quad \mbox{ for }1 \leq p \leq \infty .$$
\end{definition}
For $p=1$, this is twice the usual total variation distance. Furthermore, by Jensen's inequality, the function $p \mapsto d_{p, \pi}(\mu, \nu)$ is non-decreasing.

\begin{definition}
The $L^p$ mixing time $\tau_p(\epsilon)$ for the continuous time chain with kernel $H_t(x,y)$ and stationary distribution $\pi$ is given by 
$$ \tau_p(\epsilon) = \inf \left \{ t > 0: \sup_{x \in \mathcal{X}} d_{p,\pi}(H_t(x,\cdot), \pi) \leq \epsilon \right \}. $$
\end{definition}

Our main result bounds the $L^{\infty}$ mixing time $\tau_{\infty}(\epsilon)$, also known as the $\epsilon$-uniform mixing time. Explicitly, 
$$\tau_{\infty}(\epsilon) = \inf \left \{ t > 0 : \sup_{x,y \in \mathcal{X}} \left | \frac{H_t(x,y) - \pi(y)}{\pi(y)} \right | \leq \epsilon \right \}.$$

To estimate mixing time, we prove lower bounds on the Dirichlet form associated to the walk.
\begin{definition}
The Dirichlet form associated to $K$ is 
$$\EE_K(f,g) = \sum_{x \in \mathcal{X}} \Delta f(x) \cdot g(x) \, \pi(x) = \langle \Delta\,f,g \rangle_{\pi}$$
where $\Delta = I-K$ and $\langle \cdot,\cdot \rangle_{\pi}$ is the standard inner product for $L^2(\pi)$. 
\end{definition}
In particular, 
\begin{equation}\label{eqn:dirichlet}
\EE_K(f,f) = \frac{1}{2} \sum_{x,y \in \mathcal{X}} [f(x) - f(y)]^2 K(x,y) \pi(x).
\end{equation}
Furthermore, $\EE_K(f,f) = \EE_{K^*}(f,f)=\EE_{\frac{K+K^*}{2}}(f,f)$,
which follows from equation \eqref{eqn:dirichlet} and the identity $K^*(x,y)\pi(x)=K(y,x)\pi(y)$.
Fix $x \in \mathcal{X}$ and set $u_{x,t}(y) = h(x,y,t)$. Then recall that
\begin{equation}\label{eq:var}
\frac{d}{dt} \var(u_{x,t}) =  \sum_y \frac{d}{dt} u_{x,t}^2 d\pi 
= -2 \sum_y u_{x,t} \Delta^* u_{x,t} d\pi = -2\EE(u_{x,t}, u_{x,t}).
\end{equation}
This argument motivates the standard definition of the spectral gap
$$\lambda_1 = \inf_f \frac{\EE(f,f)}{\mbox{Var}(f)}$$
and  the well known mixing time bounds using the spectral gap:
\begin{equation} \label{eq:spec-mixing}
\tau_2(\epsilon) \leq \frac{1}{\lambda_1}\log \frac{1}{\epsilon\sqrt{\pi_*}}
\quad\textrm{and}\quad
\tau_{\infty}(1/e) \leq \frac{1}{\lambda_1} \left( 1 + \log \frac{1}{\pi_*} \right).
\end{equation}
Note that by the Courant-Fischer minmax characterization of eigenvalues, $\lambda_1$ is the second smallest eigenvalue of the symmetric operator $(\Delta + \Delta^*)/2$.

\subsection{Statement of the Main Result}

Our main result bounds the $L^{\infty}$ mixing time of a chain through eigenvalues of restricted Laplace operators.

\begin{definition}
For a non-empty subset $S \subset \mathcal{X}$, define 
$$\lambda(S) = \inf_{f \in c_0^+(S)} \frac{\EE(f, f)}{\mbox{Var}(f)}$$
where $c_0^+(S) = \{f : \mbox{supp}(f) \subset S,\,f\geq 0,\,f\neq constant\}$.
\end{definition}

In the reversible case, 
\begin{equation} \label{lambda_0}
\lambda_0(S) \leq \lambda(S) \leq \frac{1}{1-\pi(S)}\,\lambda_0(S)
\end{equation}
where $\lambda_0$ is the smallest eigenvalue of the restricted Laplacian 
$\Delta_S : c_0(S) \rightarrow c_0(S)$ with $c_0(S)=\{f : \mbox{supp}(f) \subset S\}$ and
$$ \Delta_Sf(x) = \left\{\begin{array}{cl}\Delta f(x) & x \in S \\0  & x \not \in S \end{array}\right.$$
The kernel of $\Delta_S = I - K_S$ is given explicitly by 
$$ K_S(x,y) = \left\{\begin{array}{cl} K(x,y) & x,y \in S \\0  & \mbox{otherwise} \end{array}\right.$$
By the Courant-Fischer minmax characterization of eigenvalues, \eqref{lambda_0} is equivalent to the statement:
$$
\inf_{f\in c_0(S)} \frac{\EE_{K_s}(f,f)}{\|f\|_2^2}
             \leq \inf_{f\in c_0^+(S)} \frac{\EE_K(f,f)}{\var(f)}
           \leq \frac{1}{1-\pi(S)}\,\inf_{f\in c_0(S)} \frac{\EE_{K_S}(f,f)}{\|f\|_2^2}
$$
The lower bound is due to the identity $\EE_K(f,f)=\EE_{K_S}(f,f)$ when $f\in c_0(S)$, which follows from $\Delta f(x) = \Delta_S f(x)$ when $f\in c_0(S)$ and $x\in S$. The upper bound also requires the inequality $(x-y)^2 \geq (|x|-|y|)^2$ to show that $\EE(f,f)\geq\EE(|f|,|f|)$, while Cauchy-Schwartz gives $\|f\|_1\leq\|f\|_2\sqrt{\pi(S)}$ which implies $\var(f) \geq (1-\pi(S))\|f\|_2^2$. In general, when $\pi(S)\leq 1/2$ then $\lambda(S)$ is within a factor two of the smallest eigenvalue of the symmetric operator $(\Delta_S + \Delta_S^*)/2$.
 
{\em We are interested in how $\lambda(S)$ decays as the size of $S$ increases.}

\begin{definition}
Define the {\it spectral profile} $\Lambda : [\pi_*, \infty) \rightarrow \R$ by
$$\Lambda(r) = \inf_{\pi_* \leq \pi(S) \leq r} \lambda(S)\,.$$
\end{definition}

Observe that $\Lambda(r)$ is non-increasing, and $\Lambda(r) \geq \lambda_1$. For $r \geq 1/2$, Lemma~\ref{lemma:spectral-gap} shows that $\Lambda(r)$ is within a factor two of the spectral gap $\lambda_1$.
Furthermore, by construction the walk $(K, \pi)$ satisfies the {\it Faber-Krahn} inequality
$$\lambda(S) \geq \Lambda(\pi(S)) \hspace{.5cm} \forall S \subset \mathcal{X}\,.$$

\noindent
Theorem~\ref{thm:spec-bound} is our main result:
\begin{theorem} \label{thm:spec-bound}
For $\epsilon > 0$, the $L^{\infty}$ mixing time $\tau_{\infty}(\epsilon)$ for a chain $H_t(x,y)$ satisfies
$$\tau_{\infty}(\epsilon) \leq \int_{4\pi_*}^{4/\epsilon} \frac{2dv}{v \Lambda(v)}\,.$$
\end{theorem}

In Section~\ref{sec:discrete-time}, we prove an analogous result for discrete-time walks. Since  $\Lambda(r) \geq \lambda_1$, Theorem~\ref{thm:spec-bound} shows that
$$\tau_{\infty}(1/e) \leq \int_{4\pi_*}^{4e} \frac{2\,dv}{v \Lambda(v)} 
         \leq \frac{2}{\lambda_1} \left (1 + \log \frac{1}{\pi_*} \right).$$
But since we can expect $\Lambda(r) \gg \lambda_1$ for small $r$, 
Theorem~\ref{thm:spec-bound} offers an improvement over the standard spectral gap mixing time bound \eqref{eq:spec-mixing}. 
In particular, by a discrete version of the Cheeger inequality of differential geometry, 
$$\Phi_*^2(r)/2 \leq \Lambda(r) \leq 2\Phi_*(r)$$ 
where $\Phi_*(r)$ is the (truncated) conductance profile (see Section~\ref{sec:cond}).
Consequently, by Theorem~\ref{thm:spec-bound}:
\begin{corollary} \label{cor:cond-bound}
For $\epsilon > 0$, the $L^{\infty}$ mixing time $\tau_{\infty}(\epsilon)$ for a chain $H_t(x,y)$ satisfies
$$
\tau_{\infty}(\epsilon) \leq \int_{4\pi_*}^{4/\epsilon} \frac{4dv}{v \Phi^2_*(v)}.
$$
\end{corollary}

Theorem 13 of Morris and Peres \cite{Morris} is a factor two weaker than this.

Although Theorem~\ref{thm:spec-bound} implies mixing time estimates in terms of conductance, it is reasonable to expect that for many models $\Lambda(r) \gg \Phi^2_*(r)$. In these cases, compared to Corollary~\ref{cor:cond-bound}, presently the best known conductance bound, our spectral approach leads to sharper mixing time results.
We provide below examples of such cases (see Sections~\ref{sec:examples} and \ref{sec:viscek}).

\section{Upper Bounds on Mixing Time} \label{sec:bounds}

\subsection{Spectral Profile Bounds}
In this section, we prove one of the main results, Theorem~\ref{thm:spec-bound}. The proof uses the techniques of \cite{Grig1} for estimating heat kernel decay on non-compact manifolds. The first Dirichlet eigenvalue $\lambda_0(S)$ for small sets $S$ captures the convergence behavior at the start of the walk, when the fact that the state space is finite has minimal influence. The spectral gap $\lambda_1$ governs the long-term convergence. The spectral profile $\Lambda(r)$ takes into account these two effects, since $\lambda(S) \approx \lambda_0(S)$ for $\pi(S) \leq 1/2$, and $\Lambda(r) \approx \lambda_1$ for $r \geq 1/2$.

To bound mixing times, we first lower bound $\EE(f,f)$ in terms of the spectral profile $\Lambda$, and as such Lemma~\ref{lem:Ebound} is the crucial step in the proof of Theorem~\ref{thm:spec-bound}.

We regularly use the notation that, given a function $f$, $f_+ = f \,\vee \, 0$ denotes its positive part, and $f_- = -(f \wedge 0)$ its negative part.

\begin{lemma} \label{lem:Ebound}
For every non-constant function $u:\,\mathcal{X}\mapsto \R_+$,
$$
\frac{\EE(u,u)}{\var\,u} \geq \frac 12\, \Lambda\Bigl(\frac{4(\E u)^2}{\var\,u}\Big)\,.
$$
\end{lemma}

\begin{proof}
For $c$ constant, $\EE(u,u)=\EE(u-c,u-c)$. Also, $\forall a,b\in\R:\,(a-b)^2 \geq (a_+-b_+)^2$
so $\EE(f,f)\geq \EE(f_+,f_+)$. It follows that when $0\leq c<\max u$ then
\begin{eqnarray*}
\EE(u,u) &\geq& \EE((u-c)_+,(u-c)_+) \\
 &\geq& \var((u-c)_+)\,
            \inf_{f\in c_0^+(u>c)} \frac{\EE(f,f)}{\var(f)} \\
 &\geq& \var((u-c)_+)\,
            \Lambda(\pi(u>c)) \,.
\end{eqnarray*}

Now, $\forall a,b\geq 0:\,(a-b)_+^2 \geq a^2 - 2b\,a$ and $(a-b)_+\leq a$ so
$$\var((u-c)_+) =  \E (u-c)_+^2 - (\E (u-c)_+)^2 
      \geq  \E u^2 - 2c\,\E u - (\E u)^2.$$
Let $c=\var(u)/4\E u$ and apply Markov's inequality $\pi(u>c)< (\E u)/c$, 
$$\EE(u,u) \geq (\var(u) - 2c\,\E u)\,\Lambda(\E u / c) 
   = \frac 12\,\var(u)\,\Lambda\left(\frac{4(\E u)^2}{\var\,u}\right)
$$
\end{proof}

Now we bound the $L^2$ distance of a chain from equilibrium in terms of the function $V(t): [0, \infty) \rightarrow \R$ given by
$$t = \int_{4\pi_*}^{V(t)} \frac{dv}{v \Lambda(v)}.$$
Since the integral diverges, $V(t)$ is well-defined for $t \geq 0$. 

The $L^2$ bound of Theorem~\ref{thm:L2} implies the $L^{\infty}$ bound that is our main result. To prove the $L^2$ bound, we simply apply Lemma~\ref{lem:Ebound} to the heat kernel $h(x,y,t)$.

\begin{theorem} \label{thm:L2}
For the chain $(K,\pi)$, we have 
$$\sup_{x \in \mathcal{X}} d_{2, \pi}^2(H_t(x, \cdot), \pi) \leq \frac{4}{V(t)}.$$
\end{theorem}
\begin{proof}
Given $x \in \mathcal{X}$ a value where the supremum occurs, define $u_{x,t}(y) = h(x,y,t)$ and $I_x(t) = \var(u_{x,t})$.
If $u_{x,t}=1$ then the theorem follows trivially. Otherwise, $u_{x,t}$ is non-constant
and since $\E u_{x,t} = 1$, then by \eqref{eq:var} and Lemma~\ref{lem:Ebound}
\begin{equation}
I_x'(t) = -2\EE(u_{x,t},u_{x,t}) \leq  -I_x \Lambda(4/I_x).
\end{equation}
Integrating over $[0, t]$ we have
$$\int_{I_x(0)}^{I_x(t)} \frac{dI_x}{I_x \Lambda (4/I_x)} \leq -t.$$
With the change of variable $v = 4/I_x$,
$$t \leq \int_{4/I_x(0)}^{4/I_x(t)} \frac{dv}{v \Lambda(v)}.$$
Since $I_x(0) = 1/\pi(y) - 1 < 1/\pi_*$
$$V(t) \leq \frac{4}{I_x(t)} = \frac{4}{\|h(x,\cdot, t) - 1\|_2^2}$$
and the result follows.
\end{proof}

Now we show how to transfer the $L^2$ bounds of Theorem~\ref{thm:L2} to the $L^{\infty}$ bounds of our main result. 

\noindent \textbf{Proof of Theorem~\ref{thm:spec-bound}.} 
Observe that
\begin{eqnarray}
\left| \frac{H_t(x,y)-\pi(y)}{\pi(y)}\right| &=& 
        \left| \frac{\sum_z \left(H_{t/2}(x,z)-\pi(z)\right)\left(H_{t/2}(z,y)-\pi(y)\right)}{\pi(y)} \right| \nonumber\\
    &=& \left| \sum_z \pi(z)\,\left(\frac{H_{t/2}(x,z)}{\pi(z)}-1\right)\left(\frac{H_{t/2}^*(y,z)}{\pi(z)}-1\right)\right|\nonumber \\
    &\leq& d_{2,\pi}(H_{t/2}(x,\cdot),\pi)\,d_{2,\pi}(H_{t/2}^*(y,\cdot),\pi) \label{eqn:L_infty}
\end{eqnarray}
where the inequality follows from Cauchy-Schwartz.
Since we can apply Theorem~\ref{thm:L2} to either $H_t$ or $H_t^*$, we have 
$$\sup_{x,y\in\mathcal{X}} |h(x,y,t)-1| \leq \frac{4}{V(t/2)}.$$
So $|h(x,y,t) - 1| \leq \epsilon$ for $V(t/2) \geq 4/\epsilon$, that is, for $t$ such that 
$$t/2 \geq \int_{4\pi_*}^{4/\epsilon} \frac{dv}{v \Lambda(v)}$$
proving the result. $\square$

The next result shows that any improvement in using the spectral profile $\Lambda(r)$ instead of the spectral gap $\lambda_1$ comes from looking at small sets since for $r=1/2$, already $\Lambda(r) \approx \lambda_1$.

\begin{lemma} \label{lemma:spectral-gap}
The spectral gap $\lambda_1$ and the spectral profile $\Lambda(r)$ satisfy
$$\lambda_1 \leq \Lambda(1/2) \leq 2\lambda_1.$$
\end{lemma}
\begin{proof}
The lower bound follows immediately from the definition of the spectral gap. For the upper bound, let $m$ be a median of $f$. Then using Lemma~\ref{lemma:abs},
\begin{eqnarray*}
\EE(f, f) & = & \EE(f - m, f - m) \\
& \geq & \EE((f - m)_+, (f-m)_+) + \EE((f-m)_-, (f-m)_-).
\end{eqnarray*}
Since $\pi(\{f > m\}) = \pi(\{f < m\}) \leq 1/2$, we have
$$
\EE((f - m)_+, (f-m)_+) \geq \|(f-m)_+\|_2^2 \,\lambda_0(\{f > m\}) 
$$
and
$$
\EE((f - m)_-, (f-m)_-) \geq \|(f-m)_-\|_2^2 \,\lambda_0(\{f < m\}) 
$$
Consequently, 
\begin{eqnarray*}
\EE(f,f) & \geq & \|f-m\|_2^2 \,\inf_{\pi(S)\leq 1/2} \lambda_0(S)\\
& \geq & \mbox{Var}(f) \,\frac{\Lambda(1/2)}{2} \,.
\end{eqnarray*}
The upper bound follows by minimizing over $f$.
\end{proof}

The proof required the following lemma.

\begin{lemma} \label{lemma:abs}
Given a function $f:\,\mathcal{X}\mapsto\R$ then
$$\EE(f,f) \geq \EE(f_+, f_+) + \EE(f_-, f_-) \geq \EE(|f|, |f|).$$
\end{lemma}
\begin{proof}
Given $g,h:\,\mathcal{X}\mapsto\R$ with $g,h\geq 0$ and $(\mbox{supp}\,g)\cap(\mbox{supp}\,h)=\emptyset$ then
$$
\EE(g,h) = \sum_x g(x)h(x)\pi(x) - \sum_{x,y} g(y)h(x)K(x,y)\pi(x) \leq 0
$$
because the first sum is zero and every term in the second is non-negative. 
In particular, $f_+,f_-\geq 0$ with $(\mbox{supp}\,f_+)\cap(\mbox{supp}\,f_-)=\emptyset$,
and so by linearity 
\begin{eqnarray*}
\EE(f,f) & = & \EE(f_+ - f_-, f_+ - f_-) \\
&  =  & \EE(f_+, f_+) + \EE(f_-, f_-) - \EE(f_+, f_-) - \EE(f_-, f_+) \\
&\geq & \EE(f_+, f_+) + \EE(f_-, f_-) \\
&\geq & \EE(f_+, f_+) + \EE(f_-, f_-) + \EE(f_+, f_-) + \EE(f_-, f_+) \\
&  =  & \EE(|f|, |f|).
\end{eqnarray*}
\end{proof}

\subsection{Conductance Bounds} \label{sec:cond}

In this section, we show how to use Theorem~\ref{thm:spec-bound} to recover previous bounds on mixing time in terms of the conductance profile.

\begin{definition}
For non-empty $A, B \subset \mathcal{X}$, the flow is given by 
$$Q(A,B) = \sum_{x \in A, \, y \in B} Q(x,y)$$
where $Q(x,y) = \pi(x)K(x,y)$ can be viewed as a probability measure on $\mathcal{X} \times \mathcal{X}$. The boundary of a subset is defined by 
$$\partial S = \{x \in S : \exists y \not \in S, \ K(x,y) > 0\}$$
and $|\partial S| = Q(S, S^{c})$.
\end{definition}

Observe that
$$\pi(S)=Q(S,\mathcal{X})=Q(S,S)+Q(S,S^c)$$ 
and also
$$\pi(S)=Q(\mathcal{X},S)=Q(S,S)+Q(S^c,S).$$ 
It follows that $Q(S,S^c)=Q(S^c,S)$.

Like the spectral profile $\Lambda(r)$, the conductance profile $\Phi(r)$ measures how conductance changes with the size of the set $S$.

\begin{definition}
Define the conductance profile $\Phi : [\pi_*, 1) \rightarrow \R$ by 
$$\Phi(r) = \inf_{\pi_* \leq \pi(S) \leq r} \frac{| \partial S |}{\pi(S)}$$
and the truncated conductance profile $\Phi_* : [\pi_*, 1) \rightarrow \R$ by
$$ \Phi_*(r) = \left\{\begin{array}{cc}\Phi(r) & r < 1/2 \\\Phi(1/2) & r \geq 1/2\end{array}\right.$$ 
\end{definition}
The value $\Phi(1/2)$ is often referred to as the conductance, or the isoperimetric constant, of the chain.

The next lemma is a discrete version of the ``Cheeger inequality" of differential geometry, and will let us apply Theorem~\ref{thm:spec-bound} to recover the conductance 
profile bound of Corollary~\ref{cor:cond-bound}. The proof of the lemma is similar to the proof given in \cite{FMC} of the fact that 
$$\frac{\Phi^2(1/2)}{8} \leq \lambda_1 \leq 2\Phi(1/2).$$

\begin{lemma} \label{lemma:cheeger}
For $r \in [\pi_*, 1)$, the spectral profile $\Lambda$ and the conductance profile $\Phi$ satisfy
$$\frac{\Phi^2(r)}{2} \leq \Lambda(r) \leq \frac{\Phi(r)}{1-r}.$$
\end{lemma}

\noindent
\begin{proof} It suffices to show that $\frac 12\,\Phi^2(\pi(A))\leq\lambda_0(A)\leq\frac{|\partial A|}{\pi(A)}$ for every $A\subset\mathcal{X}$. The bound then follows from \eqref{lambda_0} by minimizing over sets with $\pi(A)\leq r$.

For the upper bound,
$$\lambda_0(A)
   \leq \frac{\EE(1_A, 1_A)}{\|1_A\|_2^2} 
   = \frac{|\partial A|}{\pi(A)}.
$$

To show the lower bound, for a non-negative function $f$, define the level sets $F_t = \{ x \in \mathcal{X} : f(x) \geq t \}$ and the indicator functions $f_t = 1_{F_t}$. Then
\begin{eqnarray}
\label{eq:level} \pi(f) & = & \sum_{x \in \mathcal{X}} \left( \int_0^{\infty} f_t(x) dt \right )\pi(x) \\
& = & \int_0^{\infty} \pi(F_t) \, dt. \nonumber
\end{eqnarray}
Furthermore,
\begin{eqnarray}
\nonumber \sum_{x,y} |f(x) - f(y)| \, Q(x,y) & = & \frac{1}{2} \sum_{x,y} |f(x) - f(y)| \, [Q(x,y) + Q(y,x)] \\
& = & \sum_{f(x)> f(y)} [f(x) - f(y)] \, [Q(x,y) + Q(y,x)] \nonumber\\
& = & \sum_{f(x)> f(y)} \int_0^{\infty} 1_{\{f(y) < t \leq f(x)\}} \, [Q(x,y) + Q(y,x)] \, dt \nonumber \\
& = & \int_0^{\infty} |\partial F_t| \, dt + \int_0^{\infty} |\partial F_t^c| \, dt \nonumber \\
& = & 2 \int_0^{\infty} |\partial F_t| \, dt. \label{eq:coarea}
\end{eqnarray}
Observe that \eqref{eq:coarea} is a discrete analog of the co-area formula. For non-negative $f \in c_0(A)$, $F_t \subset A$ for $t > 0$, and so 
\begin{eqnarray*}
\sum_{x,y} |f(x) - f(y)| Q(x,y) & = & 2 \int_0^{\infty} |\partial F_t| \, dt \hspace{1cm} \mbox{by } \eqref{eq:coarea}\\
& \geq & 2 \Phi(\pi(A)) \int_0^{\infty} \pi(F_t) \, dt \\
& = & 2 \Phi(\pi(A))\pi(f) \hspace{1cm} \mbox{by } \eqref{eq:level}.
\end{eqnarray*}
Consequently, for any non-negative $f \in c_0(A)$, by the above
\begin{eqnarray*}
2 \Phi(\pi(A)) \pi(f^2) & \leq & \sum_{x,y} |f^2(x) - f^2(y)| \, Q(x,y) \\
& = & \sum_{x,y} |f(x) - f(y)|\cdot(f(x) + f(y))\, Q(x,y) \\
& \leq & \left(\sum_{x,y} (f(x) - f(y))^2 \, Q(x,y) \right)^{1/2} \\
&  & \hspace{1cm} \times \left(\sum_{x,y} (f(x) + f(y))^2 \, Q(x,y) \right)^{1/2} \\
& \leq & \left (2 \EE(f,f) \right)^{1/2} ( 4 \pi(f^2) )^{1/2}.
\end{eqnarray*}
Then
$$
\lambda_0(A) = \inf_{f\in c_0^+(A)} \frac{\EE(f,f)}{\pi(f^2)}
  \geq \frac{\Phi^2(\pi(A))}{2}\,.
$$
The infinum for $\lambda_0(A)$ occured at $f\in c_0^+(A)$ because 
for general $f\in c_0(A)$, $\EE(f,f)\geq\EE(|f|,|f|)$ and $\pi(f^2)=\pi(|f|^2)$.
\end{proof}

\begin{remark}\label{rem:cond-bound}
From the proofs of Lemma~\ref{lemma:cheeger} and
Lemma~\ref{lemma:spectral-gap}, we have that
$$
\frac{\Phi^2(1/2)}{2}
\leq \inf_{\pi(A) \leq 1/2} \lambda_0(A)
\leq \lambda_1.
$$
Consequently, when $r>1/2$ then
$\Phi^2_*(r)/2 \leq \lambda_1 \leq \Lambda(r)$,
proving the conductance profile bound of Corollary~\ref{cor:cond-bound}.
\end{remark}

\subsection{Discrete-Time Walks} \label{sec:discrete-time}

In this section we consider discrete-time chains, deriving spectral profile bounds on mixing time similar to those for continuous-time walks. For $u_{x,t}(y) = h(x,y,t)$, the rate of decay of the heat operator in the continuous setting is given by 
\begin{equation} \label{eq:deriv}
\frac{d}{dt}\var(u_{x,t}) = -2 \EE(u_{x,t}, u_{x,t}).
\end{equation}
In the discrete-time setting, set $u_{x,n}(y) = k(x,y,n) = \frac{K_n(x,y)}{\pi(y)}$.
Then, since $K^*u_{x,n} = u_{x,n+1}$ and $\E(u_{x,n}) = 1$,
\begin{eqnarray}
\label{eq:discrete-diff} \var(u_{x,n+1}) - \var(u_{x,n}) & = & \langle u_{x,n+1}, u_{x,n+1}\rangle  -   \langle u_{x,n}, u_{x,n} \rangle \\
& = & -\langle (I - KK^*) u_{x,n}, \, u_{x,n} \rangle  \nonumber \\
& = & -\EE_{KK^*}(u_{x,n}, u_{x,n}) \nonumber
\end{eqnarray}
and so it is natural to consider the multiplicative symmeterizations $KK^*$ and $K^*K$. In order to relate mixing time directly to the kernel $K$ of the original walk, we use the assumption that for $\alpha > 0$
$$K(x,x) \geq \alpha \hspace{.5cm} \forall x \in \mathcal{X}.$$

Define $\Lambda_{KK^*}$ and $V_{KK^*}$ to be the analogs of $\Lambda$ and $V$ where $\EE_K(f,f)$ is replaced by $\EE_{KK^*}(f,f)$. If $KK^*$ is reducible, then $\lambda_1^{KK^*} = 0$, and so we restrict ourselves to the irreducible case. We define $\Lambda_{K^*K}$ and $V_{K^*K}$ similarly and also assume irreducibility. The following result is a discrete-time version of Theorem~\ref{thm:L2}, and its proof is analogous.

\begin{theorem} \label{thm:discrete}
For a discrete-time chain $(K, \pi)$ with $K^*K$ and $KK^*$ irreducible
$$\sup_{x \in \mathcal{X}} d^2_{2,\pi}(K_n(x, \cdot), \pi) \leq \frac{4}{V_{KK^*}(n/2)} \quad \mbox{ and } \quad \sup_{x \in \mathcal{X}} d^2_{2,\pi}(K_n^*(x, \cdot), \pi) \leq \frac{4}{V_{K^*K}(n/2)}.$$
\end{theorem}

\noindent
\begin{proof}
The second statement follows from the first by replacing $K$ by $K^*$. For fixed $x \in \mathcal{X}$, define $u_{x,n}(y) = k(x,y,n)$ and $I_x(n) = \var(u_{x,n})$. By \eqref{eq:discrete-diff} and Lemma~\ref{lem:Ebound} 
\begin{eqnarray*}
I_x(n+1) - I_x(n) & = & -\EE_{KK^*}(u_{x,n},u_{x,n}) \\
& \leq & -\frac{1}{2} I_x(n) \Lambda_{KK^*}( 4/I_x(n)).
\end{eqnarray*}
Since both $I_x(n)$ and $\Lambda_{KK^*}(r)$ are non-increasing, the piecewise linear extension of $I_x(n)$ to $\mathbb{R_+}$ satisfies
$$I_x'(t) \leq - \frac{1}{2} I_x(t)\Lambda_{KK^*} (4/I_x(t)).$$
At integer $t$, we can take either the derivative from the right or the left.
Solving this differential equation as in Theorem~\ref{thm:L2}, we have
$$V_{KK^*}(t/2) \leq \frac{4}{I_x(t)}$$
and the result follows.
\end{proof}

\begin{corollary} \label{cor:discrete-time}
Assume that $K(x,x) \geq \alpha > 0$ for all $x \in \mathcal{X}$. Then for $\epsilon > 0$, the $L^{\infty}$ mixing time for the discrete-time chain $K_n$ satisfies
$$\tau_{\infty}(\epsilon) \leq 2\,\left\lceil \int_{4\pi_*}^{4/\epsilon} \frac{dv}{\alpha v \Lambda(v)}\right\rceil.$$
\end{corollary}

\noindent
\begin{proof}
Since $K^*(x,x) = K(x,x) \geq \alpha$, observe that
\begin{eqnarray*}
KK^*(x,y)\pi(x) & \geq & K^*(x,x)K(x,y)\pi(x) + K^*(x,y)K(y,y)\pi(x)\\
& \geq & \alpha K(x,y)\pi(x) + \alpha K(y,x)\pi(y)
\end{eqnarray*}
and so,
$$\EE_{KK^*}(f,f) \geq 2\alpha \EE_K(f,f).$$
Consequently, $\Lambda_{KK^*} \geq 2 \alpha \Lambda$, from which it follows that
\begin{eqnarray*}
\alpha t & = & \int_{4\pi_*}^{V(\alpha t)} \frac{dv}{v \Lambda(v)} \\
& \geq & 2\alpha \int_{4\pi_*}^{V(\alpha t)} \frac{dv}{v \Lambda_{KK^*}(v)}.
\end{eqnarray*}
Accordingly, $V_{KK^*}(t/2) \geq V(\alpha t)$, and similarly
$V_{K^*K}(t/2)\geq V(\alpha t)$. As in Theorem~\ref{thm:spec-bound},
\begin{eqnarray*}
|k(x,y,2n) - 1| & \leq & d_{2,\pi}(K_n(x,\cdot),\pi)\,d_{2,\pi}(K_n^*(y,\cdot),\pi) \\
& \leq & \frac{4}{V(\alpha n)}.
\end{eqnarray*}
And so, $|k(x,y,2n)-1| \leq \epsilon$ for
$$n \geq \int_{4\pi_*}^{4/\epsilon} \frac{dv}{\alpha v \Lambda(v)}.$$
\end{proof}

\noindent
{\bf Improvement for  discrete-time using rescaling}.
Given a Markov kernel $K$ let $\Lambda^K(r)$ and $\Phi^K(r)$ denote the spectral and conductance profiles, respectively. Then
\begin{eqnarray*}
\lefteqn{\Lambda^K(r) = (1-\alpha)\Lambda^{\frac{K-\alpha I}{1-\alpha}}(r)
\geq \frac{1-\alpha}{2}\,\Phi^{\frac{K-\alpha I}{1-\alpha}}(r)^2  } \\
 &=&    \frac{1-\alpha}{2}\,\left(\frac{\Phi^K(r)}{1-\alpha}\right)^2
=    \frac{\Phi^K(r)^2}{2(1-\alpha)}
\end{eqnarray*}
The appropriate discrete time version of Corollary~\ref{cor:cond-bound} is
then
\begin{corollary} \label{cor:cond-bound2}
For $\epsilon > 0$, the $L^{\infty}$ mixing time $\tau_{\infty}(\epsilon)$ for the chain $K_n$ satisfies
$$
\tau_{\infty}(\epsilon) \leq 2\left\lceil\int_{4\pi_*}^{4/\epsilon} \frac{2\,dv}{\frac{\alpha}{1-\alpha}\,v \Phi^2_*(v)} \right\rceil.
$$
\end{corollary}
In contrast, the bound of Morris and Peres \cite{Morris} is 
$$
\tau_{\infty}(\epsilon) \leq 2\left\lceil\int_{4\pi_*}^{4/\epsilon} \frac{2\,dv}{\min\left\{\frac{\alpha^2}{(1-\alpha)^2},1\right\}\,v\Phi^2_*(v)} \right\rceil
$$
which is similar, but slightly weaker when $\alpha\neq 1/2$.

\section{Lower Bounds on Mixing Time} \label{sec:lower}
In this section, we recall a result of \cite{Grig2} to show that for reversible chains the spectral profile describes well the decay behavior of the heat kernel $h_t(x,y) = H_t(x,y)/\pi(y)$. These results are based on the idea of {\it anti-Faber-Krahn inequalities}.

For reversible chains, \eqref{eqn:L_infty} implies
$$
\sup_{x,y} \frac{H_t(x,y)}{\pi(y)}-1 \leq \sup_{x} \sum_z \pi(z)\left(\frac{H_{t/2}(x,z)}{\pi(z)}-1\right)^2
  = \sup_x \frac{H_t(x,x)}{\pi(x)}-1
$$
and so
$$\sup_{x,y \in \mathcal{X}}h_t(x,y) = \sup_{x \in \mathcal{X}}h_t(x,x)\,.$$
Lemma~\ref{lemma:LB1} gives a simple lower bound on the heat kernel.

\begin{lemma}[\cite{Grig2}] \label{lemma:LB1}
For a reversible chain $(K, \pi)$ and non-empty $S \subset \mathcal{X}$,
$$\sup_{x \in \mathcal{X}}h_t(x,x) \geq \frac{\exp(-t\lambda_0(S))}{2\pi(S)}.$$
\end{lemma}
\begin{proof}
Let $\lambda_0(S) \leq \lambda_1(S) \leq \cdots \leq \lambda_{|S|-1}(S)$ be the eigenvalues of $I - K_S$. Then $K_S$ has eigenvalues $\{1 - \lambda_i(S)\}$. Since tr$(K_S^k)$ can be written as either the sum of eigenvalues, or the sum of diagonal entries, we have \begin{eqnarray*}
\sum_{i=0}^{|S|-1} (1-\lambda_i(S))^k & = & \sum_{x \in S} K_S^k(x,x) \\
& \le & \sum_{x \in S}K_k(x,x).
\end{eqnarray*}
For $k$ even, all the terms in the first sum are non-negative, and consequently 
$$(1-\lambda_0(S))^k \leq \sum_{x \in S} K_k(x,x). $$
Finally, to bound the continuous-time kernel, note that
\begin{eqnarray*}
\pi(S) \sup_{x \in \mathcal{X}} h_t(x,x) & \geq & \sum_{x \in S} h_t(x,x) \pi(x) \\
& \geq & \sum_{x \in S} e^{-t} \sum_{k=0}^{\infty} \frac{t^{2k}}{(2k)!} K_{2k}(x,x) \\
& \geq & e^{-t} \sum_{k=0}^{\infty} \frac{t^{2k}(1-\lambda_0(S))^{2k}}{(2k)!} \\
& = & e^{-t}\frac{\exp[t(1-\lambda_0(S))] + \exp[-t(1-\lambda_0(S))]}{2}
\end{eqnarray*}
from which the result follows.
\end{proof}
Theorem~\ref{thm:LB2} is a partial converse of the upper bound given in Theorem~\ref{thm:L2} under the restriction of $\delta$-regularity.

\begin{definition}
A positive, increasing function $f \in C^1(0, T)$ is $\delta$-regular if for all $0 < t < s \leq 2t < T$
$$\frac{f'(s)}{f(s)} \geq \delta \frac{f'(t)}{f(t)}.$$
\end{definition}

\begin{definition}
The walk $(K, \pi)$ satisfies the anti-Faber-Krahn inequality with function $L:[\pi_*, \infty) \to \R$ if for all $r\in [\pi_*,\infty)$, 
$$\inf_{\pi_* \leq \pi(S) \leq r} \lambda_0(S) \leq L(r).$$
\end{definition}

\begin{remark}
Observe that $(K, \pi)$ satisfies the anti-Faber-Krahn inequality with  $L(r) = \Lambda(r)$, in light of (\ref{lambda_0}).
\end{remark}

\begin{theorem}[\cite{Grig2}]\label{thm:LB2}
Let $(K, \pi)$ be a reversible Markov chain that satisfies the anti-Faber-Krahn inequality with $L : (\pi_*, \infty) \rightarrow \R$, 
and that $\gamma(t)$, defined implicitly by
$$t = \int_{\pi_*}^{\gamma(t)} \frac{dv}{v L(v)},$$
is $\delta$-regular on $(0,T)$. Then for $t \in (0,\delta T/2)$
$$\sup_{x \in \mathcal{X}} h_t(x,x) \geq \frac{1}{2 \gamma(2t/\delta)}.$$
\end{theorem}

\begin{proof}
Fix $t \in (0, \delta T/2)$ and set $r = \gamma(t/\delta)$. By the anti-Faber-Krahn inequality, there exists $S \subset \mathcal{X}$ with $\pi(S) \leq r$ and $\lambda_0(S) \leq L(r)$. Consequently, by Lemma~\ref{lemma:LB1},
$$\sup_{x \in \mathcal{X}}h_t(x,x) \geq \frac{\exp(-t\lambda_0(S))}{2\pi(S)} \geq \frac{\exp(-tL(r))}{2r}.$$
So, $\sup_x h_t(x,x) \geq \exp(-C_t)$ for $C_t = \log 2r + t L(r)$. Since $L(\gamma(s)) = (\log \gamma)' (s)$
$$C_t = \log 2\gamma(t/\delta) + t (\log \gamma)'(t/\delta).$$
By the mean value theorem, there exists $\theta \in (t/\delta, 2t/\delta)$ such that 
$$(\log \gamma)'(\theta) = \frac{\log\gamma(2t/\delta) - \log \gamma(t/\delta)}{t/\delta}.$$
By $\delta$-regularity 
$$(\log \gamma)'(\theta) \geq \delta(\log \gamma)'(t/\delta)$$ 
and so $C_t \leq \log[2 \gamma(2t/\delta)]$, showing the result.
\end{proof}

\section{Applications} \label{sec:apps}

The following lemma, while hardly surprising, is often effective in reducing computation in specific examples. In particular, it is used in computing the spectral profile of the random walk on the $n$-cycle in the present section.

\begin{lemma} \label{lemma:connected}
Let $S = S_1 \cup \dots \cup S_k$ be a decomposition of $S$ into connected components. Then 
$$\lambda(S)= \min_{S_i} \{\lambda(S_i)\}.$$
\end{lemma}

\begin{proof}
Clearly $\lambda(S) \leq \min_{S_i} \{\lambda(S_i)\}$, and we need only show the reverse inequality. For a function $f \geq 0$, define $f_{S_i} = 1_{S_i}f$. Then
$$
\var(f) = \var\left(\sum_{S_i} f_{S_i} \right)
 = \sum_{S_i} \E f_{S_i}^2 - \left( \sum_{S_i} \E f_{S_i} \right)^2 
\leq  \sum_{S_i} \var(f_{S_i}).
$$
Consequently,
\begin{eqnarray*}
\lambda(S) & = & \inf_{f \in c_0^+(S)} \frac{\EE(f,f)}{\var(f)} \\
& = & \inf_{f \in c_0^+(S)} \frac{\sum_{S_i} \EE(f_{S_i}, f_{S_i})}{\var(f)}\\
& \geq & \inf_{f \in c_0^+(S)} \frac{\sum_{S_i} \lambda(S_i) \var(f_{S_i})}{\var(f)}
\end{eqnarray*}
and the result follows.
\end{proof}

\subsection{First Examples}\label{sec:examples}
\subsubsection{The Complete Graph}

Consider the continuous-time walk on the complete graph in the $n$-point space $\Omega=\{x_1, \dots, x_n\}$ with kernel $K(x_i, x_j) = 1/n$ $\forall i,j$.
To find the eigenvalues of the restricted operator $K_S: c_0(S) \mapsto c_0(S)$, we consider functions $f: \{x_1, \dots, x_{|S|}\} \mapsto \R$. Since 
$$K_Sf(x_j) = \frac{1}{n} \sum_{i=1}^{|S|} f(x_i) = \bar{f} \hspace{1cm} 1 \leq j \leq |S|$$
$f$ is an eigenfunction of $K_S$ with corresponding eigenvalue $\lambda$ if and only if $\lambda f(x_j) = \bar{f}$ for $1 \leq j \leq |S|$. If $\lambda  \not = 0$, then this implies that $f$ is constant with eigenvalue $\lambda = |S|/n$.
So, the smallest eigenvalue of $I - K_S$ satisfies $\lambda_0(S) = 1 - |S|/n$, and the second smallest eigenvalue of $I-K$ satisfies $\lambda_1=1$.
Since 
$$\lambda_1 \leq \lambda(S) \leq \frac{\lambda_0(S)}{1-\pi(S)}$$
$\lambda(S) = 1$ and accordingly $\Lambda(r) \equiv 1$. 

Theorem~\ref{thm:spec-bound} then shows that for the complete graph $\tau_{\infty}(\epsilon) \leq 2 \log (n/\epsilon)$. Since the distribution of the chain at any time $t \geq 0$ is given explicitly by 
$$H_t(x_i, x_j) = e^{-t} \delta_{x_i}(x_j) + \frac{(1 - e^{-t})}{n}$$
we see that $\tau_{\infty}(\epsilon) = \log[(n-1)/\epsilon]$, and so our estimate is off by a factor of 2.

\subsubsection{The $n$-Cycle} \label{sec:ncycle}
To bound the spectral profile $\Lambda(r)$ for simple random walk on the $n$-cycle,
recall (see Lemma~\ref{lemma:connected})  that it is sufficient to restrict our attention to connected subsets. 
Now consider simple random walk on the $n$-cycle $\Omega = \{x_0, \dots x_{n-1}\}$ given by kernel $K(x_i, x_j) = 1/2$ if $j = i \pm 1 \pmod{n}$ and zero otherwise. By Lemma~\ref{lemma:connected}, to find $\lambda_0(S)$ we need only consider connected subsets $S \subset \Omega$. For $S$ such that $\pi(S) < 1$, $I-K_S$ corresponds to the tridiagonal Toeplitz matrix
with 1's along the diagonal and 1/2's along the upper and lower off-diagonals
(and 0's everywhere else).
In this case, the least eigenvalue is given explicitly by 
$$\lambda_0(S) = 1 - \cos\left(\frac{\pi}{|S|+1}\right).$$
Since the spectral gap satisfies 
$\lambda_1 = 1 - \cos ({2\pi}/{n})$, we have 
$\Lambda(r) \approx 1/(rn)^2$ for $1/n \leq r \leq 1$. Theorem~\ref{thm:spec-bound} then shows the correct $O(n^2)$ mixing time bound.

\subsection{Log-Sobolev and Nash Inequalities}
Logarithmic Sobolev and Nash inequalities are among the strongest tools available to study $L^2$ convergence rates of finite Markov chains. Log-Sobolev inequalities were introduced by Gross \cite{Gross2,Gross} to study Markov semigroups in infinite dimensional settings, and developed in the discrete setting by Diaconis and Saloff-Coste \cite{LogSob}. Nash inequalities were originally formulated to study the decay of the heat kernel in certain parabolic equations (see \cite{Nash}). Building on ideas in \cite{Carlen,SC-Nash, SC-Nash2}, Diaconis and Saloff-Coste \cite{DSNash} show how to apply Nash's argument to finite Markov chains. In this section we show that both log-Sobolev and Nash inequalities yield bounds on the spectral profile $\Lambda(r)$, leading to new proofs of previous mixing time estimates in terms of these inequalities.

\begin{definition}
The log-Sobolev constant $\rho$ is given by
$$\rho = \inf_{\ent_{\pi} f^2\neq 0} \frac{\EE(f,f)}{\ent_{\pi} f^2}$$
where the entropy $\displaystyle \ent_{\pi}(f^2) = \sum_{x \in \mathcal{X}}f^2(x)\log\left(f^2(x)/\|f\|_2^2\right)\pi(x).$
\end{definition}

\begin{lemma}[Log-Sobolev] \label{lemma:logsob}
The spectral profile $\Lambda(r)$ and log-Sobolev constant $\rho$ satisfy
$$
\Lambda(r) \geq \rho\,\frac{\log(1/r)}{1-r}\,.
$$
\end{lemma}

\begin{proof}
By definition
$$
\Lambda(r) = \inf_{\pi(S)\leq r} \inf_{f\in c_0^+(S)} \frac{\EE(f,f)}{\var_{\pi}(f)}
  \geq \rho\,\inf_{\pi(S)\leq r} \inf_{f\in c_0^+(S)} \frac{\ent_{\pi}(f^2)}{\var_{\pi}(f)}
$$
The lemma will follow if for every set $S\subset\mathcal{X}$
$$
\inf_{f\in c_0^+(S)} \frac{\ent_{\pi}(f^2)}{\var_{\pi}(f)} 
\geq \frac{\log \frac{1}{\pi(S)}}{1-\pi(S)}\,.
$$
Define a probability measure $\pi'(x)=\frac{\pi(x)}{\pi(S)}$ if $x\in S$ and $\pi'(x)=0$ otherwise.  Then
$$
\inf_{f\in c_0^+(S)} \frac{\ent_{\pi}(f^2)}{\var_{\pi}(f)} 
= \inf_{f\in c_0^+(S)}
  \frac{\ent_{\pi'}(f^2) + \log \frac{1}{\pi(S)}\,E_{\pi'} f^2}{E_{\pi'} f^2 - \pi(S)\,(E_{\pi'}f)^2}
$$
Rearranging the terms, it suffices to show that
$$
\inf_{f\in c_0^+(S)} \frac{\ent_{\pi'}(f^2)}{\var_{\pi'}(f)} \geq \frac{\pi(S)\log \frac{1}{\pi(S)}}{1-\pi(S)}\,.
$$
However, since $\pi(S)\in(0,1)$ then $\pi(S)\frac{\log(1/\pi(S))}{1-\pi(S)}\leq 1$ and so it suffices that for every probability measure and $f\geq 0$ that $\ent(f^2)/\var(f)\geq 1$. This is true, as observed in \cite{LO00.1} and recalled in Remark 6.7 of 
\cite{BT03.1}.
\end{proof}

\noindent
The bound $\lambda_0(A)\geq\rho\log(1/\pi(A))$ can be shown similarly, but without need for the result
of \cite{LO00.1}.
Like log-Sobolev inequalities, Nash inequalities also yield bounds on the spectral profile:

\begin{lemma}[Nash Inequality] \label{lemma:nash}
Given a Nash inequality
$$
\|f\|_2^{2+1/D} 
  \leq C\,\left[\EE(f,f) + \frac{1}{T}\,\|f\|_2^2\right]\,\|f\|_1^{1/D}
$$
which holds for every function $f:\,\mathcal{X}\mapsto\R$ and some constants
$C,\,D,\,T\in\R_+$, then
$$
\Lambda(r) \geq \frac{1}{C\,r^{1/2D}} - \frac{1}{T}.
$$
\end{lemma}

\begin{proof}
The Nash inequality can be rewritten as 
$$
\frac{\EE(f,f)}{\|f\|_2^2} \geq 
    \frac{1}{C}\,\left(\frac{\|f\|_2}{\|f\|_1}\right)^{1/D} - 
    \frac{1}{T}
$$
Then, 
\begin{eqnarray*}
\lambda_0(A) = \inf_{f\in c_0(A)} \frac{\EE(f,f)}{\|f\|_2^2} 
  &\geq& \inf_{f\in c_0(A)} \frac{1}{C}\,\left(\frac{\|f\|_2}{\|f\|_1}\right)^{1/D} - \frac{1}{T} \\
  &\geq& \frac{1}{C\,\pi(A)^{1/2D}} - \frac{1}{T}\,.
\end{eqnarray*}
The final inequality was due to Cauchy-Schwartz:
$\|f\|_1 \leq \|f\|_2\,\sqrt{\pi(\mbox{supp}\,f)}$.
The lemma follows by minimizing over $\pi(A)\leq r$.
\end{proof}

Although the spectral profile $\Lambda(r)$ is controlled by the spectral gap $\lambda_1$ for $r \geq 1/2$, Nash inequalities tend to be better for $r$ close to $0$, and log-Sobolev inequalities for intermediate $r$. Combining Lemmas \ref{lemma:logsob} and \ref{lemma:nash}, we get the following bounds on mixing time:
\begin{corollary}
Given the spectral gap $\lambda_1$ and the log-Sobolev constant $\rho$ and/or 
a Nash inequality with $DC\geq T$, $D\geq 1$ and $\pi_*\leq 1/4e$, the $L^{\infty}$ mixing time for the continuous-time Markov chain 
with $\epsilon\leq 8$ satisfies
\begin{eqnarray*}
\tau_{\infty}(\epsilon) &\leq &
\frac{2}{\rho}\,\log \log \frac{1}{4\pi_*} + \frac{2}{\lambda_1}\,\log \frac{8}{\epsilon} \\
\tau_{\infty}(\epsilon) &\leq&
4T + \frac{2}{\lambda_1}\,\left(2D\log \frac{2DC}{T} + \log \frac{4}{\epsilon}\right) \\
\tau_{\infty}(\epsilon) &\leq&
4T + \frac{2}{\rho}\,\log \log \left(\frac{2DC}{T}\right)^{2D} 
   + \frac{2}{\lambda_1}\,\log \frac{8}{\epsilon}
\end{eqnarray*}
\end{corollary}


\begin{proof}
For the first upper bound use the log-Sobolev bound $\Lambda(r)\geq\rho\log(1/r)$
when $r<1/2$ and the spectral gap bound when $r\geq 1/2$. Simple integration gives the result.

For the second upper bound use the Nash bound when $r\leq (T/2DC)^{2D}$ and spectral gap bound for the remainder. Then
\begin{eqnarray*}
\tau_{\infty}(\epsilon) 
&\leq& 
\int_{4\pi_*}^{(T/2DC)^{2D}} 
   \frac{2dr}{r\,\frac{1}{C\,r^{1/2D}}\,\left(1-\frac{C\,r^{1/2D}}{T}\right)} 
+ \int_{(T/2DC)^{2D}}^{4/\epsilon} \frac{2dr}{r\,\lambda_1} \\
&\leq&
4T + \frac{2}{\lambda_1}\,\log \frac{4/\epsilon}{(T/2DC)^{2D}} 
\end{eqnarray*}
where the second inequality used the bound
$1-\frac{C\,r^{1/2D}}{T} \geq 1-\frac{1}{2D} \geq 1/2$
before integrating. Simplification gives the result.

For the mixed bound use the Nash bound when $r\leq (T/2DC)^{2D}$, the log-Sobolev bound for
$(T/2DC)^{2D}\leq r< 1/2$ and the spectral gap bound when $r\geq 1/2$.
\end{proof}

Similar discrete time bounds follow from Corollary~\ref{cor:discrete-time}. When $\forall x:\,K(x,x)\geq\alpha$ then these bounds are roughly a factor $\alpha^{-1}$ larger than the continuous time case.

These bounds compare well with previous results shown through different methods.
For instance, Aldous and Fill \cite{AldousBook} combine results of Diaconis and Saloff-Coste \cite{LogSob,DSNash} to show the continuous time bound
$$
\tau_{\infty}(\epsilon) \leq
2T + \frac{1}{\rho}\,\log \log\left(\frac{DC}{T}\right)^D
  + \frac{1}{\lambda_1}\,\left(4+\log(1/\epsilon)\right)
$$
whenever $DC\geq T$.

\subsubsection{Walks with Moderate Growth} \label{sec:modgrowth}
In this section, we describe how estimates on the volume growth of a walk give estimates on the spectral profile $\Lambda(v)$. The treatment given here is analogous to the method of Nash inequalities described in \cite{DSNash}. 

Define the Cayley graph of $(K, \pi)$ to be the undirected graph on the state space $\mathcal{X}$ with edge set $E = \{(x,y) : \pi(x)K(x,y) + \pi(y)K(y,x) > 0 \}$. Let $d(x,y)$ be the usual graph distance, and denote the closed ball of radius $r$ around $x$ by $B(x,r) = \{z : d(x,z) \leq r \}$. The volume of $B(x,r)$ is given by $V(x,r) = \sum_{z \in B(x,r)} \pi(z)$.

\begin{definition} For $A,d \geq 1$, the finite Markov chain $(K,\pi)$ has $(A,d)$-moderate growth if 
\begin{equation} \label{eq:modgrowth}
V(x,r) \geq \frac{1}{A}\left(\frac{r+1}{\gamma}\right)^d \hspace{.5cm}\forall x \in \mathcal{X}, \  0 \leq r \leq \gamma
\end{equation}
where $\gamma$ is the diameter of the graph.
\end{definition}

For any $f$ and $r \geq 0$, set
$$f_r(x) = \frac{1}{V(x,r)} \sum_{y \in B(x,r)}f(y) \pi(y).$$
\begin{definition}
The finite Markov chain $(K, \pi)$ satisfies a local Poincar{\'e} inequality with constant $a$ if for all $f$ and $r \geq 0$
\begin{equation} \label{eq:poincare}
\|f - f_r\|_2^2 \leq a r^2 \EE(f,f).
\end{equation}
\end{definition}

Under assumptions \eqref{eq:modgrowth} and \eqref{eq:poincare}, Diaconis and Saloff-Coste \cite{DSNash} derive the Nash inequality
\begin{equation} \label{eq:mod-growth-nash}
\|f\|_2^{2 + 4/d} \leq C\left[ \EE(f,f) + \frac{1}{a\gamma^2} \|f\|_2^2 \right]\|f\|_1^{4/d}
\end{equation}
where $C = (1+1/d)^2(1+d)^{2/d}A^{2/d}a\gamma^2$. By Lemma~\ref{lemma:nash}, this immediately implies the lower bound on the spectral profile
$$\Lambda(v) \geq \left(\frac{d^2}{(d+1)^{2 + 2/d}A^{2/d}v^{2/d}} - 1\right) \frac{1}{a\gamma^2}.$$ 

Theorem~\ref{thm:moderate-growth} below shows how to bound $\Lambda(r)$ in 
terms of a local Poincar{\'e} inequality and the volume growth function
$$V_*(r) = \inf_x V(x,r).$$
The proof is similar to the derivation of Nash inequalities for walks with 
moderate growth shown in \cite{DSNash}.

\begin{theorem} \label{thm:moderate-growth}
Let $(K, \pi)$ be a finite Markov chain that satisfies the local 
Poincar{\'e} inequality with constant $a$. For $v \leq 1/2$, the spectral 
profile satisfies
$$\Lambda(v) \geq \frac{1}{4aW^2(2v)}$$
where $W(v) = \inf \{r : V_*(r) \geq v \}$.
\end{theorem}

\begin{proof}
Fix $S \subset \mathcal{X}$ with $\pi(S) \leq 1/2$ and $f \in c_0(S)$. It 
is sufficient to show that
$$\frac{\EE(f,f)}{\|f\|_2^2} \geq \frac{1}{4aW^2(2\pi(S))}.$$
First observe that
\begin{eqnarray*}
\|f\|_2^2 & = & \langle f-f_r, f \rangle + \langle f_r, f \rangle \\
& \leq & \|f-f_r\|_2 \cdot \|f\|_2 + \langle f_r, f \rangle.
\end{eqnarray*}
Now,
\begin{eqnarray*}
\langle f_r, f \rangle  & = & \sum_x \left ( \frac{1}{V(x,r)} \sum_{y \in B(x,r)} f(y) 
\pi(y) \right) f(x) \pi(x) \label{eq:r-average}\\
& \leq & \frac{1}{V_*(r)} \|f\|_1^2 \\
& \leq & \frac{\pi(S)}{V_*(r)} \|f\|_2^2.
\end{eqnarray*}
Consequently, by the local Poincar{\'e} inequality,
$$\|f\|_2^2 \leq \sqrt{a}r\EE(f,f)^{1/2}\|f\|_2 + 
\frac{\pi(S)}{V_*(r)}\|f\|_2^2.$$
Dividing by $\|f\|_2^2$ and choosing $r = W(2\pi(S))$ we have
$$1 \leq \sqrt{a}r\frac{\EE(f,f)^{1/2}}{\|f\|_2} + 1/2$$
and the result follows.
\end{proof}


\begin{corollary} \label{cor:modmixing}
Let $(K, \pi)$ be a finite Markov chain that satisfies $(A,d)$-moderate 
growth and the local Poincar{\'e} inequality with constant $a$. Then the 
$L^{\infty}$ mixing time satisfies
$$\tau_{\infty}(\epsilon) \leq  C(a,A,d,\epsilon) \gamma^2$$
where $\gamma$ is the diameter of the graph and $C(a,A,d,\epsilon)$ is a 
constant depending only on $a$, $A$, $d$ and $\epsilon$.
\end{corollary}
\begin{proof}
By the moderate growth assumption, $W(v) \leq \gamma(Av)^{1/d}$. And so 
for $v \leq 1/2$
$$\Lambda(v) \geq \frac{1}{4aW^2(2v)} \geq \frac{1}{8aA^{1/d}\gamma^2 
v^{2/d}}.$$
For $v \geq 1/2$, note that $\Lambda(v) \geq \lambda_1 \geq 
\Lambda(1/2)/2$. The result now follows immediately from 
Theorem~\ref{thm:spec-bound}.
\end{proof}

In Theorem 3.1 of \cite{ModGrowth} Diaconis and Saloff-Coste show that for walks on groups
with $(A,d)$-moderate growth, local Poincar\'e inequality with constant $a$,
and $\gamma\geq A4^{d+1}$
$$\tau_{\infty}(1/e) \geq \frac{\gamma^2}{4^{2d+1}A^2}\,.$$
It follows that $\tau_{\infty}(1/e)=\Theta(\gamma^2)$, and Corollary \ref{cor:modmixing} was of the correct order $\gamma^2$.

For instance, consider the example of simple random walk on the $n$-cycle discussed in Section~\ref{sec:ncycle}. For this walk $V(x_i, r) = (1+2\lfloor r \rfloor)/n$, and so it satisfies the moderate growth criterion \eqref{eq:modgrowth} with $A = 6$, $d=1$ and diameter $\gamma = \lfloor n/2 \rfloor$. Moreover, it is shown in \cite{DSNash} that every group walk satisfies the local Poincar{\'e} inequality
$$\|f - f_r\|_2^2 \leq 2|S|r^2\EE(f,f)$$
where $S$ is a symmetric generating set for the walk. Consequently, Corollary~\ref{cor:modmixing} shows that the walk on the $n$-cycle mixes in $O(n^2)$ time. For several additional examples of walks with moderate growth, see \cite{ModGrowth, DSNash}.

\subsection{The Viscek Graphs}\label{sec:viscek}
For a random walk $(K, \pi)$ consider its Cayley graph defined in Section~\ref{sec:modgrowth}. Define the minimum volume of a disk of radius $r$ by $V_*(r) = \inf_x\{V(x,r)\}$. Here we first use a result of \cite{Barlow} that shows that the spectral profile $\Lambda(r)$ can be bounded in terms of the volume growth $V_*(r)$ alone
(see Lemma~\ref{lemma:vol}). We then apply this technique to analyze walks on  the fractal-like Viscek family of finite graphs.  

\begin{lemma}[\cite{Barlow}] \label{lemma:vol}
Let $Q_* = \inf_{x \sim y}[\pi(x)K(x,y) + \pi(y)K(y,x)]$ and 
$$w(r) = \inf\{k : V_*(k) > r\}.$$
Then
$$\lambda(A) \geq \frac{Q_*}{4\pi(A)w(\pi(A))}.$$
\end{lemma}

\begin{proof}
Fix $f \in c_0(A)$ normalized so that $\|f\|_{\infty} = 1$. Then,
$$\|f\|_2^2 = \sum_x |f(x)|^2 \pi(x) \leq \pi(A).$$
Let $x_0$ be a point such that $|f(x_0)| = 1$ and let $k = \max\{ l \in \mathbb{N} : B(x_0, l) \subset A \}$. Then there is a sequence of points $x_0, x_1, \dots, x_{k+1}$ with $x_i \sim x_{i+1}$, $x_0, \dots, x_k \in A$ and $x_{k+1} \not \in A$. So,
\begin{eqnarray*}
\EE(f,f) & = & \frac{1}{2} \sum_{x,y} |f(x) - f(y)|^2 \pi(x)K(x,y) \\
& \geq & \frac{1}{2}\sum_{i = 0}^k|f(x_{i+1}) - f(x_i)|^2 [\pi(x_i)K(x_i,x_{i+1}) + \pi(x_{i+1})K(x_{i+1}, x_i)]\\
& \geq & \frac{Q_*}{2(k+2)}\left (\sum_{i = 0}^k|f(x_{i+1}) - f(x_i)| \right) ^2 \\
& = & \frac{Q_*}{2(k+2)} |f(x_{k+1}) - f(x_0)|^2 \\
& = & \frac{Q_*}{2(k+2)}.\\
\end{eqnarray*}
Consequently, 
$$\lambda_0(A) = \inf_{f \in c_0(A)} \frac{\EE(f,f)}{\|f\|_2^2} \geq \frac{Q_*}{2(k+2)\pi(A)}.$$
To finish the proof, observe that $\pi(A) \geq V(x_0, k) \geq V_*(k)$, and so $w(\pi(A)) \geq k+1 \geq (k+2)/2$.
\end{proof}

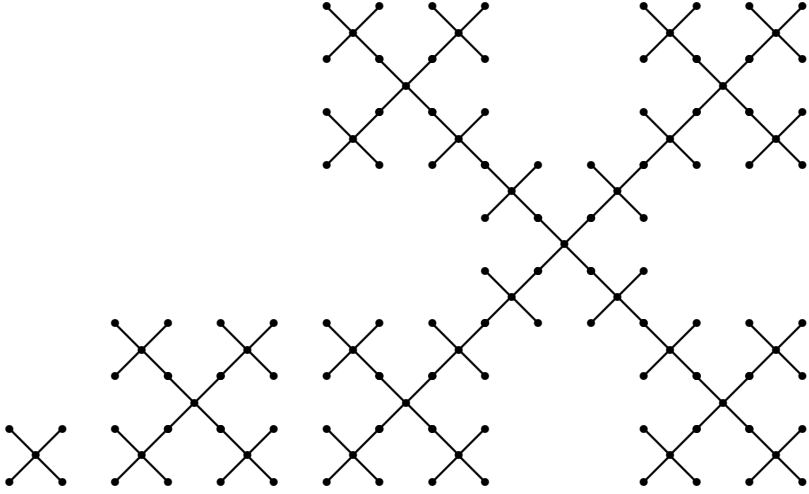
\begin{figure}
\begin{picture}(300,200)
\thicklines
\put(0,20){\line(1,-1){20}} 
\put(0,0){\line(1,1){20}} 
\multiput(0,0)(10,10){3}{\circle*{3}}
\multiput(0,20)(20,-20){2}{\circle*{3}}

\multiput(40,0)(20,20){3}{
\put(0,20){\line(1,-1){20}} 
\put(0,0){\line(1,1){20}} 
\multiput(0,0)(10,10){3}{\circle*{3}}
\multiput(0,20)(20,-20){2}{\circle*{3}}
}
\multiput(40,40)(40,-40){2}{
\put(0,20){\line(1,-1){20}} 
\put(0,0){\line(1,1){20}} 
\multiput(0,0)(10,10){3}{\circle*{3}}
\multiput(0,20)(20,-20){2}{\circle*{3}}
}

\multiput(80,0)(60,60){3}{
\multiput(40,0)(20,20){3}{
\put(0,20){\line(1,-1){20}} 
\put(0,0){\line(1,1){20}} 
\multiput(0,0)(10,10){3}{\circle*{3}}
\multiput(0,20)(20,-20){2}{\circle*{3}}
}
\multiput(40,40)(40,-40){2}{
\put(0,20){\line(1,-1){20}} 
\put(0,0){\line(1,1){20}} 
\multiput(0,0)(10,10){3}{\circle*{3}}
\multiput(0,20)(20,-20){2}{\circle*{3}}
}}

\multiput(80,120)(120,-120){2}{
\multiput(40,0)(20,20){3}{
\put(0,20){\line(1,-1){20}} 
\put(0,0){\line(1,1){20}} 
\multiput(0,0)(10,10){3}{\circle*{3}}
\multiput(0,20)(20,-20){2}{\circle*{3}}
}
\multiput(40,40)(40,-40){2}{
\put(0,20){\line(1,-1){20}} 
\put(0,0){\line(1,1){20}} 
\multiput(0,0)(10,10){3}{\circle*{3}}
\multiput(0,20)(20,-20){2}{\circle*{3}}
}}
\end{picture}
\caption{The first three generations $\mathcal{V}_4(0)$, $\mathcal{V}_4(1)$ and $\mathcal{V}_4(2)$ of a Viscek graph with $N=4$.}
\label{fig:viscek}
\end{figure}

The Viscek graphs are a two parameter family of finite trees that are inductively defined as follows. Fix the parameter $N \geq 2$, and define $\mathcal{V}_N(0)$ to be the star graph on $N+1$ vertices (i.e. a central vertex surrounded by $N$ vertices). Given $\mathcal{V}_N(n-1)$ choose $N$ vertices $x_1, \dots, x_N$ such that $d(x_i, x_j) = \mbox{diam}(\mathcal{V}_N(n-1))$ for $i\neq j$. Construct $\mathcal{V}_N(n)$ by taking $N+1$ copies $\{\mathcal{V}_N^i(n-1)\}_{i=0}^N$ of the $(n-1)^{th}$ generation graph, and for $1 \leq i \leq N$ identifying $x_i^0 \in \mathcal{V}_N^0(n-1)$ with $x_i^i \in \mathcal{V}_N^i(n-1)$. Observe that a different choice of vertices $x_1, \dots, x_N$ leads to an isomorphic construction. For $N=2$, $\mathcal{V}_2(n)$ is a path for each $n$. Figure~\ref{fig:viscek} illustrates the first three generations of a Viscek graph for $N=4$.

The following lemma bounds the spectral profile and mixing time for simple random walk on $\mathcal{V}_N(n)$. The proof is analogous to the volume growth computation for the infinite Viscek graph $\mathcal{V}_N(\infty) = \lim_{n\rightarrow \infty} \mathcal{V}_N(n)$ given in \cite{Barlow} and recalled in \cite{SC-FK}. 

\begin{lemma}\label{lemma:viscek}
For $N \geq 2$, $r \leq 1$ the spectral profile $\Lambda(r)$ for simple random walk on $\mathcal{V}_N(n)$ satisfies
$$\frac{a(N)}{\gamma^{d+1} r^{1+1/d}} \leq \Lambda(r) \leq \frac{A(N)}{\gamma^{d+1} r^{1+1/d}} \hspace{1cm} d=\log_3(N+1)$$
where $\gamma = \mbox{diam}(\mathcal{V}_N(n)) = 2\cdot 3^n$ and the constants $a, A > 0$ depend only on $N$. 

In particular, there exist constants $b,B > 0$ depending only on $N$ such that the mixing time for the continuous-time walk satisfies
$$b(N)\gamma^{d+1} \leq \tau_1(1/e) \leq \tau_{\infty}(1/e) \leq B(N)\gamma^{d+1}.$$
\end{lemma}

Observe that since the conductance profile for $\mathcal{V}_N(n)$ satisfies 
$$\Phi(r) \approx \frac{1}{|E_N(n)|r} \approx \frac{1}{\gamma^d r},$$ 
using the conductance profile bound of Corollary~\ref{cor:cond-bound} results in the upper bound $\tau_{\infty}(1/e) \preceq \gamma^{2d}$ which overestimates the mixing time for $N \geq 3$.

\begin{proof}
We first show that the mixing time bound follows from the spectral profile estimate. The upper bound is a direct consequence of Theorem~\ref{thm:spec-bound}. Recall that for an ergodic chain, the spectral gap $\lambda_1$ and $L^1$ mixing time are related by $1/\lambda_1 \leq \tau_1(1/e)$ (see e.g. \cite{FMC}). Since $\Lambda(r) \geq \lambda_1$, the lower bound is immediate.

To estimate the spectral profile, first note that the number of edges $|E_N(n)| = N(N+1)^n$. Since $\mathcal{V}_N(n)$ is a tree, $|\mathcal{V}_N(n)| = N(N+1)^n + 1$. Furthermore, $\mbox{diam}(\mathcal{V}_N(n)) = 2\cdot 3^n$.

For $0 \leq k \leq n$, define a $k$-block to be a subgraph of $\mathcal{V}_N(n)$ isomorphic to the $k^{th}$ generation graph $\mathcal{V}_N(k)$. Fix $x \in \mathcal{V}_N(n)$ and $3 \leq r \leq \mbox{diam}(\mathcal{V}_N(n))$. Then there is a unique integer $m$ such that $3^{m+1} \leq r < 3^{m+2}$. Moreover, the vertex $x$ is contained in some $m$ block $B$. Since diam$(B) = 2\cdot 3^m$, $B(x,r) \supseteq B$. Consequently, 
$$|B(x,r)| \geq |B| = N(N+1)^m + 1$$
and since $\pi_* = 1/(2|E_N(n)|)$
$$V_*(r) \geq \frac{N(N+1)^m + 1}{2N(N+1)^n} \succeq \left(\frac{r}{\gamma}\right)^d$$
where $d = \log_3(N+1)$ and the notation $a \preceq b$ indicates that there is some constant $c(N) > 0$ depending only on $N$ such that $a \leq c(N)b$. Thus, using the notation of Lemma~\ref{lemma:vol}, $w(s) \preceq \gamma s^{1/d}$. Since $Q_* = 1/|E_N(n)| \succeq 1/\gamma^d$, Lemma~\ref{lemma:vol} gives the lower bound on the spectral profile.

For the upper bound we construct test functions $f_m$ supported on $m$-blocks. Given an $m$-block $A \subset \mathcal{V}_N(n)$, choose vertices $x_1, \dots, x_N$ such that $d(x_i, x_j) = \mbox{diam}(A)$ for $i \neq j$, and call the shortest paths between these vertices diagonals. These diagonals meet in a unique point $o$ at the center of the $m$-block. Define the function $f_m \in c_0(A)$ as follows: Along diagonals, $f_m$ varies linearly with $f_m(o) = 1$ and $f_m(x_i) = 0$. Since $d(o,x_i) = \mbox{diam}(A)/2 = 3^m$, along diagonals the function is given explicitly by $f_m(x) = 1 - d(o,x)/3^m$. For a point $x$ off of the diagonals, let $f_m(x) = f_m(x')$ where $x'$ is the closest point to $x$ that lies on a diagonal. (See Figure~\ref{fig:viscek2} for a graphical representation of $f_m$). Now, since $K(x,y)\pi(x) = 1/(2|E_N(n)|)$ for $x \sim y$
\begin{eqnarray*}
\EE(f_m, f_m) & = & \frac{1}{2} \sum_{x,y} |f_m(x) - f_m(y)|^2 K(x,y)\pi(x) \\
& = & 3^{-2m} \cdot \frac{N3^m}{2|E_N(n)|}\\
& \approx & \frac{1}{\gamma^d 3^m}.
\end{eqnarray*}

\begin{figure}
\begin{picture}(180,200)
\thicklines

\multiput(0,0)(60,60){3}{
\multiput(0,0)(20,20){3}{
\put(0,20){\line(1,-1){20}} 
\put(0,0){\line(1,1){20}}
}
\multiput(0,40)(40,-40){2}{
\put(0,20){\line(1,-1){20}} 
\put(0,0){\line(1,1){20}} 
}}

\put(0,0){\circle*{2}} \put(10,10){\circle*{3}} \put(20,20){\circle*{4}}
\put(30,30){\circle*{5}} \put(40,40){\circle*{6}} \put(50,50){\circle*{7}}
\put(60,60){\circle*{8}} \put(70,70){\circle*{9}} \put(80,80){\circle*{10}}
\put(90,90){\circle*{11}}

\put(180,180){\circle*{2}} \put(170,170){\circle*{3}} \put(160,160){\circle*{4}}
\put(150,150){\circle*{5}} \put(140,140){\circle*{6}} \put(130,130){\circle*{7}}
\put(120,120){\circle*{8}} \put(110,110){\circle*{9}} \put(100,100){\circle*{10}}

\multiput(0,20)(10,-10){3}{\circle*{3}} \multiput(160, 180)(10,-10){3}{\circle*{3}}
\multiput(0,60)(10,-10){7}{\circle*{5}} \multiput(120,180)(10,-10){7}{\circle*{5}}
\multiput(0,40)(120,120){2}{\multiput(0,0)(40,-40){2}{\multiput(0,0)(10,10){3}{\circle*{5}}}}

\multiput(40,60)(80,80){2}{\multiput(0,0)(10,-10){3}{\circle*{7}}}
\multiput(60,80)(40,40){2}{\multiput(0,0)(10,-10){3}{\circle*{9}}}

\multiput(0,120)(120,-120){2}{
\multiput(0,0)(20,20){3}{
\put(0,20){\line(1,-1){20}} 
\put(0,0){\line(1,1){20}} 
}
\multiput(0,40)(40,-40){2}{
\put(0,20){\line(1,-1){20}} 
\put(0,0){\line(1,1){20}} 
}}

\put(0,180){\circle*{2}} \put(10,170){\circle*{3}} \put(20,160){\circle*{4}}
\put(30,150){\circle*{5}} \put(40,140){\circle*{6}} \put(50,130){\circle*{7}}
\put(60,120){\circle*{8}} \put(70,110){\circle*{9}} \put(80,100){\circle*{10}}

\put(180,0){\circle*{2}} \put(170,10){\circle*{3}} \put(160,20){\circle*{4}}
\put(150,30){\circle*{5}} \put(140,40){\circle*{6}} \put(130,50){\circle*{7}}
\put(120,60){\circle*{8}} \put(110,70){\circle*{9}} \put(100,80){\circle*{10}}

\multiput(0,160)(10,10){3}{\circle*{3}} \multiput(160, 0)(10,10){3}{\circle*{3}}
\multiput(0,120)(10,10){7}{\circle*{5}} \multiput(120,0)(10,10){7}{\circle*{5}}
\multiput(0,140)(120,-120){2}{\multiput(0,0)(40,40){2}{\multiput(0,0)(10,-10){3}{\circle*{5}}}}

\multiput(40,120)(80,-80){2}{\multiput(0,0)(10,10){3}{\circle*{7}}}
\multiput(60,100)(40,-40){2}{\multiput(0,0)(10,10){3}{\circle*{9}}}

\end{picture}

\caption{A graphical representation of the test function $f_2$ supported on a 2-block.}
\label{fig:viscek2}
\end{figure}
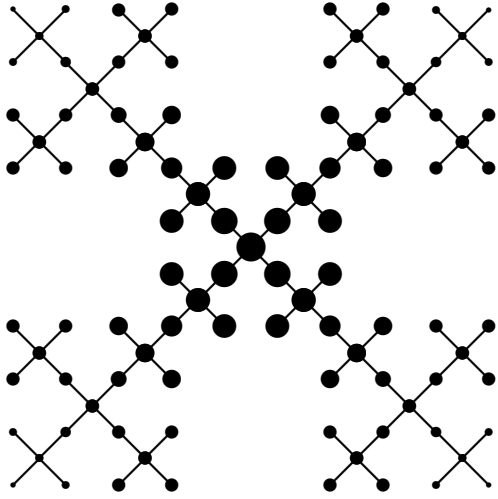

Define the central $m-1$ block of $A$ to be $A' = \{x \in A : d(o,x) \leq 3^{m-1}\}$. Since $f_m(x) \geq 2/3$ on $A'$,
$$\|f_m\|_2^2 \geq \frac{4}{9} \pi(A') \approx \frac{(N+1)^m}{\gamma^d}.$$
It is sufficient to prove the upper bound on $\Lambda(r)$ for $1/(N+1)^{n-2} < r \leq 1/2$. For these $r$ take
$$m(r) = \left \lfloor \frac{\log r(N+1)^{n-2}}{\log N+1} \right \rfloor \leq n.$$
Then $(N+1)^{m(r)} \leq r(N+1)^{n-2}$ and so for an $m(r)$-block $K$, $\pi(K) \leq r$. Consequently, for $r$ in this range
\begin{eqnarray*}
\frac{\EE(f_m, f_m)}{\var(f_m)} & \leq & \frac{2\EE(f_m, f_m)}{\|f_m\|_2^2} \\
& \preceq & \frac{1}{[3(N+1)]^m} 
\end{eqnarray*} 
Finally, since $(N+1)^m \succeq r\gamma^d$
\begin{eqnarray*}
\frac{\EE(f_m,f_m)}{\var(f_m)} & \preceq & (N+1)^{-m\frac{\log 3(N+1)}{\log N+1}} \\
& = & (N+1)^{-m(1 + 1/d)} \\
& \preceq & \frac{1}{\gamma^{d+1}r^{1+1/d}}
\end{eqnarray*}
and the upper bound on $\Lambda(r)$ follows.
\end{proof}


\subsection{A delicate example}\label{sec:delicate_eg}

Consider simple random walk on the product group $\mathbb{Z}_n \times 
\mathbb{Z}_{n^2}$. For this model, it is not hard to see that $\gamma = 
\Theta(n^2)$ and that the volume satisfies

$$V_*(r) \asymp \left\{ \begin{array}{ll} (r+1)^2/n^3 & 0 \leq r \leq n \\ 
r/n^2 & n \leq r \leq n^2 \end{array} \right. .$$
Taking $r = 0$ in \eqref{eq:modgrowth} shows that walks with moderate 
growth must have
$$\frac{1}{n^3} \geq \frac{1}{A}\left(\frac{1}{n^2}\right)^d.$$
Consequently, $\mathbb{Z}_n \times \mathbb{Z}_{n^2}$ is of moderate growth 
with $d=3/2$ and furthermore, this is the optimal choice of $d$ (assuming 
$A$ and $d$ are constant). Corollary~\ref{cor:modmixing} gives the correct 
$\gamma^2 = n^4$ mixing time, but gives the underestimate
$$ \Lambda(v) \geq \frac{C(a,A)}{\gamma^2 v^{4/3}}$$
for the spectral profile. The problem is that the moderate growth 
criterion alone is not sufficient to identify the two different 
scales of volume growth present in this example: For $r \ll 1/n$ the 
space appears $2$-dimensional, while for $r \gg 1/n$ it looks $1$-dimensional. However, we can apply 
Theorem~\ref{thm:moderate-growth} to directly take into account volume 
estimates, leading to sharp bounds on both the spectral profile and the rate 
of decay of $d_{\infty,\pi}(H_t, \pi)$.

\begin{lemma} \label{lemma:prod-walk}
For $a,b \geq 2$, the walk on $G = \mathbb{Z}_a \times \mathbb{Z}_b$ with generating set $\{(\pm 1, 0), (0,\pm 1)\}$ has spectral profile satisfying
$$\Lambda(v) \asymp \left\{ \begin{array}{ll} 1/vab & 1/ab \leq v \leq a/b 
\\ 1/v^2 b^2 & a/b \leq v \leq 1 \\ 1/b^2 & 1 \leq v \end{array} \right. .$$
In particular,
$$d_{\infty,\pi}(H_t, \pi) \asymp \left\{ \begin{array}{ll} ab / (t+1) & 0 \leq t \leq a^2 \\ 
b/t^{1/2} & a^2 \leq t \leq b^2 \end{array} \right.$$
and there are constants $c_1, c_2 > 0$ such that for $t \geq b^2$
$$e^{-c_1t/b^2} \preceq d_{\infty, \pi}(H_t,\pi) \preceq e^{-c_2t/b^2}.$$
\end{lemma}

\begin{proof}
Without loss of generality, assume $a \leq b$. Then diam$(G) = \Theta(b)$, and the volume function is given by 
$$V_*(r) \asymp \left\{ \begin{array}{ll} (r+1)^2/(ab) & 0 \leq r \leq a
 \\ r/b & a \leq r \leq \mbox{diam}(G) \end{array} \right. .$$
Consequently, $W(v) = \inf \{r : V_*(r) \geq v \}$ satisfies
$$W(v) \asymp \left\{ \begin{array}{ll} v^{1/2}(ab)^{1/2} & 1/ab \leq v \leq a/b 
\\ vb & a/b \leq v \leq 1 \end{array} \right. .$$
Since this group walk is driven by a constant number of generators, the 
chain satisifes a local Poincar{\'e} inequality with constant independent of $a$ and $b$. For $v \leq 1/2$, the lower bound on the $\Lambda(v)$ now follows from Theorem~\ref{thm:moderate-growth}. By Lemma~\ref{lemma:spectral-gap}, $\lambda_1 \geq \Lambda(1/2)/2$. For $v \geq 1/2$, the lower bound on $\Lambda(r)$ then follows from the fact $\Lambda(r) \geq \lambda_1$. 

For the upper bound, for $m \geq 1$  define linear functions $f_m \in 
c_0(B(0,m))$ by
$$f_m(x) = 1 - \frac{d_G(0,x)}{m}.$$
Then,
\begin{eqnarray*}
\EE(f_m, f_m) & = & \frac{1}{2|G|}\sum_{x,y}[f_m(x) - f_m(y)]^2K(x,y) \\
& \leq & \frac{V_*(m)}{2m^2}.
\end{eqnarray*}
Since $G$ is volume doubling (i.e. $V_*(2m) \approx V_*(m)$)
$$\|f_m\|_2^2 \geq \left(\frac{1}{2}\right)^2 V_*(m/2) \succeq V_*(m)$$
and consequently, for $V_*(m) \leq 1/2$
$$\frac{\EE(f_m,f_m)}{\var(f_m)} \leq 2 \frac{\EE(f_m,f_m)}{\|f_m\|_2^2} 
\preceq \frac{1}{m^2}.$$
Since it is sufficient to consider only $v \leq 1/2$, the upper bound on 
$\Lambda(v)$ follows by taking $m = W(v) -1$.

Following the notation of Theorem~\ref{thm:LB2}, define $\gamma(t)$ implicitly through  
$$L(r) = \left\{ \begin{array}{ll} C/rab & 1/ab \leq r \leq 
a/b \\ C/r^2b^2 & a/b \leq r \leq 1 \end{array} \right. $$
and observe that for $C$ sufficiently large, the walk satisfies the anti-Faber-Krahn inequality with $L(r)$ by the upper bound on $\Lambda(r)$. Then $\gamma(t)$ is given explicitly as 
$$\gamma(t) = \left\{ \begin{array}{ll} (Ct + 1)/ab & 0 \leq t \leq 
(a^2 - 1)/C \\ (2Ct - a^2 + 2)^{1/2}/b & (a^2 -1)/C \leq  t \leq (b^2 + a^2 -2)/(2C) \end{array} \right. .$$
Furthermore, 
$$\frac{\gamma'(t)}{\gamma(t)} = \left\{ \begin{array}{ll} C/(Ct + 1) & 0 < t \leq 
(a^2 - 1)/C \\ C/(2Ct - a^2 + 2) & (a^2 -1)/C \leq  t < (b^2 + a^2 -2)/(2C) \end{array} \right. .$$
So, $\gamma(t)$ is $\delta$-regular on $(0,T)$ with $\delta = 1/6$ and $T = (b^2 + a^2 -2)/(2C)$. For $t \leq cb^2$ and $c$ sufficiently small, the lower bound on convergence now follows from Theorem~\ref{thm:LB2}. 

Let $\{\lambda_i\}$ be the eigenvalues of $I-K$ with corresponding real orthonormal eigenfunctions $\{\psi_i\}$. Since $h_t(x,y) = H_t 1_y(x)$, writing $1_y(\cdot)$ in this $L^2$ basis we have
$$h_t(x,y) = \sum_{i = 0}^{n-1} e^{-t\lambda_i}\psi_i(x) \psi_i(y) = 1 +  \sum_{i = 1}^{n-1} e^{-t\lambda_i}\psi_i(x) \psi_i(y).$$
In particular,
$$\sup_x h_t(x,x) - 1 \geq \sup_x e^{-t\lambda_1} \psi_i^2(x) \geq e^{-t\lambda_1}$$
since the eigenfunctions are normalized in $L^2$. Since $\lambda_1 \approx \Lambda(1/2) \approx 1/b^2$, the lower bound on convergence rate follows.

Now we show the upper bound. By the lower bound on $\Lambda(r)$, using the notation of Theorem~\ref{thm:L2},
$$V(t) \succeq \left\{ \begin{array}{ll} t/ab & 0 \leq t \leq a^2 \\ 
t^{1/2}/b & a^2 \leq  t \leq b^2 \\ e^{ct/b^2} & b^2 \leq t \end{array} \right. .$$
Consequently, 
$$\sup_{x \in \mathcal{X}} d_{2, \pi}^2 (H_t(x,\cdot), \pi) \preceq \left\{ \begin{array}{ll} ab/t & 0 \leq t \leq a^2 \\ 
b/t^{1/2} & a^2 \leq  t \leq b^2 \\ e^{-ct/b^2} & b^2 \leq t\end{array} \right. .$$ 
Since the walk is reversible, the upper bound on $d_{\infty, \pi}$ then follows from the argument of Theorem~\ref{thm:spec-bound}  
\end{proof}

For $a \ll b$, Lemma~\ref{lemma:prod-walk} captures the fact that decay is fast at the start of the walk and then slows down. More specifically,
$$d_{\infty,\pi}(H_t, \pi) \asymp 1/V_*(t^{1/2}) \hspace{5mm} t \leq \mbox{diam}^2(G).$$ 
As shown in \cite{DSC-Harnack}, this relationship holds in general for random walks on groups with volume doubling. While we considered only the simple case of $\mathbb{Z}_a \times \mathbb{Z}_{b}$, the technique applies well to more general $k$-fold products.

\section*{Acknowledgments}
We thank Laurent Saloff-Coste for initiating this project by pointing out that the method of Faber-Krahn inequalities developed in \cite{Grig1} also applies in the setting of finite Markov chains. His PhD student Sharad Goel thanks Prof. Saloff-Coste for his plentiful guidance, advice and encouragement during this work.

\end{document}